\documentclass[10pt,journal]{IEEEtran}
\usepackage{amssymb}
\usepackage{verbatim}
\usepackage[final]{graphicx}
\usepackage[cmex10]{amsmath}
\usepackage{cite}
\usepackage{multirow}
\usepackage{array}
\usepackage{framed}
\usepackage{epstopdf}
\usepackage{caption,subcaption}
\usepackage{color}

\def\A{{\mathbf A}}
\def\x{{\mathbf x}}
\def\X{{\mathbf X}}

\def\s{{\mathbf s}}
\def\z{{\mathbf z}}

\def\Tr{\mathrm{Trace}}

\def\C{{\mathbb C}}

\begin{document}

\title{Feasible Point Pursuit and Successive Approximation of Non-convex QCQPs}
\author{Omar Mehanna, \emph{Student Member, IEEE}, Kejun Huang, \emph{Student Member, IEEE}, Balasubramanian Gopalakrishnan, \emph{Student Member, IEEE}, Aritra Konar, \emph{Student Member, IEEE} and Nicholas D. Sidiropoulos$^\dagger$,~\emph{Fellow, IEEE}
\thanks{Submitted to {\it IEEE Signal Processing Letters}, September 11, 2014.
Supported in part by NSF ECCS-1231504, NSF AST-1247885, NSF IIS-1247632.}
\thanks{$^\dagger$ Contact author, phone: (612) 625-1242, Fax: (612) 625-4583.}}

\maketitle
\begin{abstract}
Quadratically constrained quadratic programs (QCQPs) have a wide range of applications in signal processing and wireless communications. Non-convex QCQPs are NP-hard in general. Existing approaches relax the
non-convexity using semi-definite relaxation (SDR) or linearize the non-convex part and solve the resulting convex problem. However, these techniques are seldom successful in even obtaining a feasible solution when the QCQP matrices are indefinite. In this paper, a new feasible point pursuit successive convex approximation (FPP-SCA) algorithm is proposed for non-convex QCQPs. FPP-SCA linearizes the non-convex parts of the problem as conventional SCA does, but adds
slack variables to sustain feasibility, and a penalty to ensure slacks are sparingly used. When FPP-SCA is successful in identifying a feasible point of the non-convex QCQP, convergence to a Karush-Kuhn-Tucker (KKT) point is thereafter ensured. Simulations show the effectiveness of our proposed algorithm in obtaining feasible and near-optimal solutions, significantly outperforming existing approaches.
\end{abstract}

\section{Introduction}
\label{Sec:Intro}
Quadratically constrained quadratic programs (QCQPs) are an important class of optimization problems that have a wide spectrum of applications ranging from transmit beamforming in wireless networks, to portfolio risk management in financial engineering \cite{Luo10,Pal14}. A QCQP can be expressed as
\begin{equation}\label{QCQP}
\addtolength{\fboxsep}{3pt} \boxed{ \begin{split}
   \large{\bf{\Pi_{1}}}\: \min_{\x \in \mathbb{C}^n } &  \quad \x^H \A_0 \x  \\
      \mathrm{s.t.}    &    \quad  \x^H \A_m \x \leq c_m , \quad m=1,\ldots,M  \end{split}}
\end{equation}
where ${\bf A}_0 \succeq 0$, i.e., positive semi-definite, and ${\bf A}_m \in \C^{n\times n}$ are Hermitian matrices for all $m$ $\in$ $\left\{1,\ldots,M\right\}$. In the special case when ${\bf A}_m \succeq 0$ for all $m$ $\in$ $\left\{1,\ldots,M\right\}$, the QCQP $\large{\bf{\Pi_{1}}}$ becomes a convex optimization problem which can be efficiently solved to optimality using interior point methods \cite{Boyd}. For general indefinite ${\bf A}_m$, this problem is non-convex and NP-hard \cite{BoydLec}, except for special cases, such as when $M \leq 3$ \cite{Pat98,YeZhang03,HuaPal:2010}.

Several methods have been proposed to approximate non-convex QCQPs, including (a) the (prevailing) {\em semi-definite relaxation} (SDR) approach \cite{Boyd}; (b) the {\em reformulation linearization technique} (RLT) \cite{Anstreicher12,SherAdam98}; and (c) {\em successive convex approximation} (SCA) \cite{Beck10,Marks78,Tran14,Scutari14}. RLT consists of a reformulation step and a linearization step. The reformulation step creates redundant nonlinear constraints involving pairwise product combinations of the individual scalar variables, by multiplying different constraint pairs. The linearization step then substitutes a continuous variable for each distinct product of variables. The resulting convex optimization problem is solved to obtain an approximate solution of the non-convex problem. The main issue with RLT is that the solution of the linear program is seldom feasible for the non-convex problem. Furthermore the size of the linear program approximation is much larger than the original problem, thereby making it computationally involved.

The SCA approach is a more general scheme to deal with non-convex problems, and its application to non-convex QCQPs is sometimes called the convex-concave approach \cite{BoydLec}: each quadratic term is separated into convex and concave parts, and the latter is replaced by a convex (usually linear) approximation around a feasible point. The resulting convex problem is solved to obtain the next iterate, which also serves as the approximation point for the next iteration. Scutari {\it et al.} recently proposed a parallel and distributed SCA framework to obtain stationary points for non-convex optimization problems \cite{Scutari14}. The algorithm starts from an initial feasible point, the non-convex constraints are approximated by a strictly convex function around this point, and the resulting convex problem is solved to obtain the next iterate. This procedure is repeated until convergence to a stationary point.

The drawback with approaches (b) and (c) is that they need a feasible point as initialization, which is difficult to obtain in general. Constraint approximation about a feasible point yields a nonempty set that contains at least the given point, whereas constraint approximation about an infeasible point tends to yield an empty set, even if the original problem is feasible. Existing convergence results for SCA depend on a feasible initialization, e.g., \cite{Beck10}.

The most popular among the above is the SDR approach \cite{Boyd,Luo10}, where the original problem is reformulated by introducing $\X:=\x \x^H$ and solving the semi-definite program (SDP)  $\large{\bf{\Pi_{1r}}}$ obtained after relaxing the rank-1 constraint.
 \begin{equation}\label{SDR}
 \begin{split}
    \large{\bf{\Pi_{1r}}} \: \min_{\X \in \mathbb{C}^{n\times n} } &  \quad \Tr({\bf A_0} \X)   \\
      \mathrm{s.t.}    &    \quad  \Tr(\A_m \X) \leq c_m , ~ m=1,\ldots,M \\
						&    \quad \X \succeq 0 \end{split}
\end{equation}
Because of the rank relaxation, the solution to $\large{\bf{\Pi_{1r}}}$ gives a lower bound on the optimal value of the cost function of $\large{\bf{\Pi_{1}}}$. Note that the SDR $\large{\bf{\Pi_{1r}}}$ is the Lagrange bi-dual of $\large{\bf{\Pi_{1}}}$. If the solution $\X^*$ to problem $\large{\bf{\Pi_{1r}}}$ is rank-1, the optimal solution $\x^*$ to $\large{\bf{\Pi_{1}}}$ is the principal eigenvector of $\X^*$, scaled by the square root of the maximum eigenvalue of $\X^*$; otherwise, randomization techniques are used \cite{Luo10}. If the matrices $\{\A_m\}_{m=1}^M$ are all negative semidefinite, then any randomly generated point can be scaled up to satisfy the constraints of the QCQP $\large{\bf{\Pi_{1}}}$; finding a feasible solution using randomization is easy in this case, and the challenge is to find one that is close to optimum, see \cite{Nikos06,Luo10}.

In the general setting where $\{\A_m\}_{m=1}^M$ are all indefinite, or when one deals with two-sided positive semidefinite constraints such as in \cite{Phan09},  SDR with randomization often fails to find a point that satisfies the constraints in $\large{\bf{\Pi_{1}}}$. That is why it is important to develop an alternate approach (instead of SDR followed by randomization) that has a high probability of finding a feasible solution to the NP-hard QCQP $\large{\bf{\Pi_{1}}}$, when one exists.

In this paper we propose an iterative algorithm for obtaining good feasible solutions to general QCQPs, where we approximate the feasible region through a linear restriction of the non-convex parts of the constraints. In order to guarantee feasibility of the modified problem, slack variables are added, and a penalty is used to ensure that slacks are sparingly used. The solution of the resulting optimization problem is then used to compute a new linearization, and the procedure is repeated until convergence. The proposed {\em Feasible Point Pursuit - Successive Convex Approximation} (FPP-SCA) algorithm differs from the conventional SCA approach \cite{Beck10,Marks78,Tran14} in that the latter requires the starting point to be in the feasible region of the original problem. Finding a feasible point is easy in some cases, such as when all the constraints $\x^H \A_m \x \leq c_m$ involve negative semi-definite $\A_m$, as considered in \cite{Tran14}. For general QCQPs, however, finding an initial feasible point is hard. The performance of the proposed algorithm is compared with the conventional SDR followed by randomization, and simulations show that the proposed algorithm attains a feasible solution for a much larger percentage of problem instances. Furthermore, the feasible solution obtained using the proposed algorithm is much closer to the SDR lower bound than SDR followed by randomization.

\section{The FPP-SCA approach}
Problem $\large{\bf{\Pi_{1}}}$ may or may not be feasible, and establishing (in)feasibility is generally NP-hard. When infeasible, one may instead seek a compromise that minimizes constraint violations in some sense - this is common in engineering applications. In order to account for potential infeasibility, consider adding slack variables ${\bf s} \in \mathbb{R}^M$ and a slack penalty to $\large{\bf{\Pi_{1}}}$ 
\begin{equation}\label{QCQP_slack}
\addtolength{\fboxsep}{3pt} \boxed{ \begin{split}
   \large{\bf{\Pi_{2}}}\:
   \min_{\x \in \mathbb{C}^n, \s \in \mathbb{R}^M } &
   		\quad \x^H \A_0 \x + \lambda\|{\bf s}\|  \\
      \mathrm{s.t.}    &    \quad  \x^H \A_m \x \leq c_m + s_m , \\
            & \quad s_m \geq 0, \quad m=1,\ldots,M,
        \end{split}}
\end{equation}
where $\lambda$ trades off the original objective function and the slack penalty term, and $\|\cdot\|$ can be any vector norm. Problem $\large{\bf{\Pi_{2}}}$ is always feasible, and if $(\x_o,\s_o)$ is an optimal solution of $\large{\bf{\Pi_{2}}}$ and it so happens that $\s_o = {\bf 0}$, then $\x_o$ is an optimal solution of $\large{\bf{\Pi_{1}}}$; else using the $l_1$ norm of $\s$ in $\large{\bf{\Pi_{2}}}$ (which reduces to the sum of the slacks, due to the non-negativity constraints) promotes sparsity in terms of constraint violations. The difficulty though is that problem $\large{\bf{\Pi_{2}}}$ remains non-convex and NP-hard in general.

\noindent {\bf Successive convex approximation (SCA).} Using eigen-decomposition, the matrix $\A_m$ can be expressed as $\A_m = \A_m^{(+)} + \A_m^{(-)}$, where  $\A_m^{(+)} \succeq 0$ and  $\A_m^{(-)} \preceq 0$ (negative semi-definite). For any $\z,\x \in \C^{n \times 1}$,  $(\x-\z)^H \A_m^{(-)} (\x- \z) \leq 0$. Expanding the left-hand side of the inequality, we obtain
\begin{equation}\label{linapprx1}
\x^H \A_m^{(-)} \x \leq 2 \text{Re}\left\{\z^H \A_m^{(-)} \x\right\} - \z^H \A_m^{(-)} \z .
\end{equation}
Therefore, using the linear {\em restriction} \eqref{linapprx1} around the point $\z$, we may replace the $m$-th (non-convex) constraint of $\large{\bf \Pi_2}$ with the convex constraint
\begin{equation}\label{apprx}
\x^H \A_m^+ \x + 2 \text{Re}\left\{\z^H \A_m^{(-)} \x\right\} \leq c_m + \z^H \A_m^{(-)} \z + s_m .
\end{equation}
This leads us to propose the following algorithm.
\newline

\noindent \textbf{{\sf Feasible Point Pursuit Successive Convex Approximation (FPP-SCA) Algorithm}}

\noindent \textbf{Initialization:} Set $k=0$ and randomly generate an initial point $\z_0$.

\noindent \textbf{Repeat}
\begin{enumerate}
\item Solve
\begin{equation}\label{SOCP}
\hspace{-10pt}\addtolength{\fboxsep}{3pt} \boxed{ \begin{split}
    \large{\bf{\Pi_{3}}}\; \min_{\x ,\s} &
    		\quad \x^H \A_0 \x + \lambda \sum_{m=1}^M s_m  \\
     \mathrm{s.t.}  ~ & \quad \x^H \A_m^{(+)} \x +
     			2\text{Re}\left\{\z_k^H \A_m^{(-)}\x\right\} \\
     			& \quad \quad \quad
     				\leq c_m + \z_k^H \A_m^{(-)}\z_k + s_m\\
     			      &  \quad s_m \geq 0,\quad m=1,\ldots,M \end{split}}
\end{equation}
\item Let $\x_k$ denote the optimal $\x$ obtained by solving $\large{\bf{\Pi_{3}}}$ at the $k$-th iteration, and set $\z_{k+1}=\x_k$.
\item Set $k=k+1$.
\end{enumerate}
\noindent \textbf{until convergence.}
\newline

\noindent Some important remarks and claims are in order.

\noindent $\bullet$ We first relaxed the constraints in $\large{\bf{\Pi_{1}}}$ by adding slacks, then tightened the relaxed constraints via partial linear restriction of their non-convex parts. We could instead first tighten the original constraints (risking turning a feasible original problem into an infeasible one) then relax by adding slacks to make the restriction feasible -- the net result turns out being $\large{\bf{\Pi_{3}}}$ in both cases, and it is always feasible. 

\noindent $\bullet$ FPP-SCA yields a non-increasing cost sequence, i.e., the optimal cost of $\large{\bf{\Pi_{3}}}$ is non-increasing in $k$. This is because the cost function is independent of $k$, and the solution of the $k$-th iteration is also feasible for the $(k+1)$-th iteration. To see this, note that $\z_{k+1} = \x_k$ is the optimal solution of $\large{\bf \Pi_3}$ at the $k$-th iteration, so it satisfies the restriction $\z_{k+1}^H \A_m^+ \z_{k+1} + 2 \text{Re}\left\{\z_{k}^H \A_m^{(-)} \z_{k+1}\right\} \leq c_m + \z_{k}^H \A_m^{(-)} \z_{k} + s_m$, and therefore {\em a fortiori} also the non-convex quadratic constraint $\z_{k+1}^H \A_m^+ \z_{k+1} + \z_{k+1}^H \A_m^- \z_{k+1} \leq c_m + s_m$. Looking at the corresponding constraint at the next iteration, $\x^H \A_m^+ \x + 2 \text{Re}\left\{\z_{k+1}^H \A_m^{(-)} \x\right\} \leq c_m + \z_{k+1}^H \A_m^{(-)} \z_{k+1} + s_m$, plugging in $\x=\z_{k+1}$ we obtain $\z_{k+1}^H \A_m^+ \z_{k+1} + \z_{k+1}^H \A_m^- \z_{k+1} \leq c_m + s_m$, i.e., feasibility of $(\x_k=\z_{k+1},\s)$ for the $k$-th iteration of $\large{\bf \Pi_3}$ implies feasibility of $(\z_{k+1},\s)$ for the same $\s$ for the $(k+1)$-th iteration of $\large{\bf \Pi_3}$.

\noindent $\bullet$ Problem $\large{\bf{\Pi_{3}}}$ is convex and can be easily formulated as a second-order cone program (SOCP). Worst-case complexity of solving $\large{\bf{\Pi_{3}}}$ is $\mathcal{O}\left([n+M]^{3.5}\right)$ $\ll$ $\mathcal{O}(n^{6.5})$ for SDR \cite{Boyd96}, and FPP-SCA usually takes just a few iterations to converge.

\noindent $\bullet$ FPP-SCA can be run using different starting points $\z_0$, and the best solution can be taken. Simulations suggest that SDR can provide a good initialization for FPP-SCA, if the extra complexity of SDR is acceptable. Otherwise random initialization(s) can be used.

\noindent $\bullet$ We propose using $\lambda \gg 1$ to force the slack variables toward zero, thereby pushing the iterates towards the feasible region of $\large{\bf{\Pi_{1}}}$ when this is non-empty. Higher $\lambda$ also helps ensure that if a feasible point of $\large{\bf{\Pi_{1}}}$ is found, subsequent iterates will remain in the feasible region of $\large{\bf{\Pi_{1}}}$, although there are no analytical guarantees for this. Note that the main advantage of FPP-SCA over conventional SCA is the ability to find a feasible point with high probability. Once the slacks are all 0, one can simply switch to conventional SCA without the slacks. Our simulations show that the subsequent iterates are almost identical between these two schemes (with or without the slack variables).

\noindent $\bullet$ If FPP-SCA converges, it converges to a KKT point for problem $\large{\bf{\Pi_{2}}}$, according to Beck {\it et al} \cite{Beck10}. If the converged slack variables turn out being all zero, then it is easy to show that the remaining variables satisfy the KKT conditions for the original problem $\large{\bf{\Pi_{1}}}$.

\begin{figure}[t]
\centering
\begin{subfigure}[t]{0.2\textwidth}
\includegraphics[width=\textwidth]{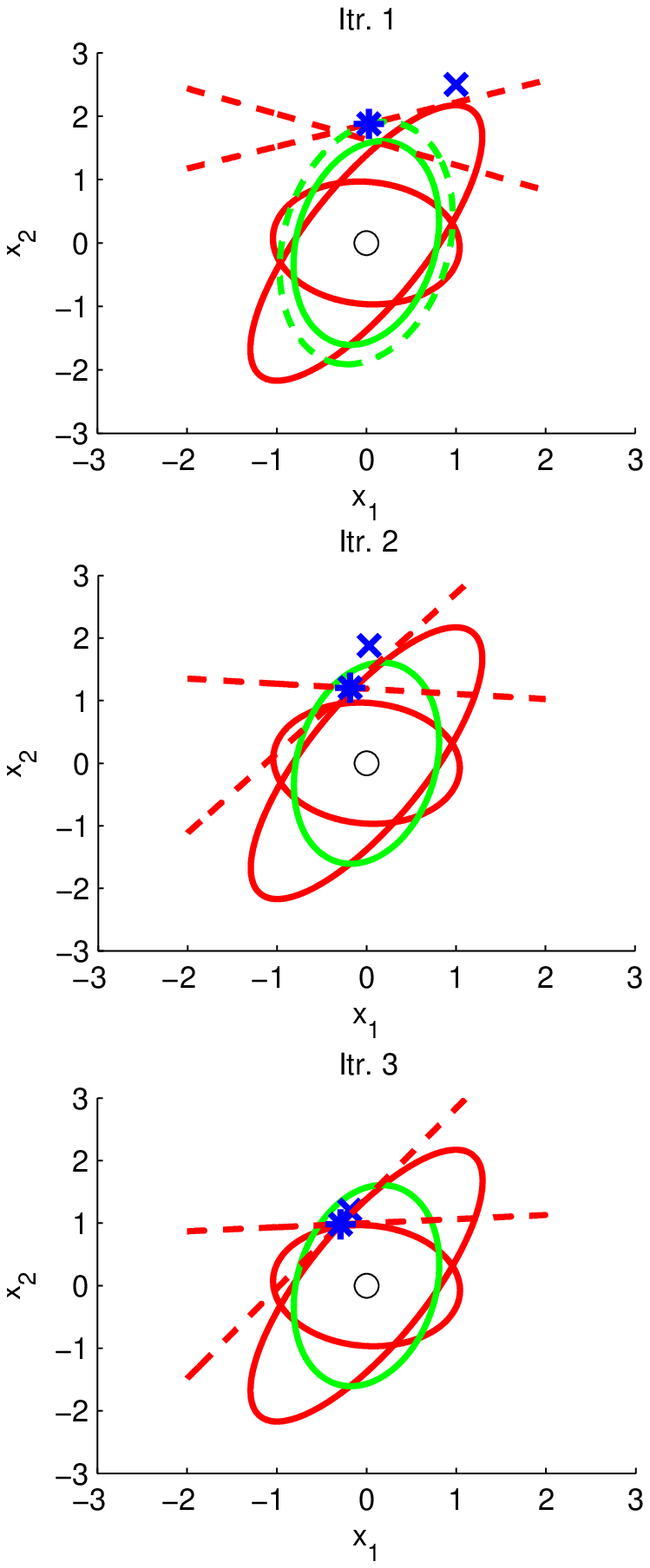}
\caption{Successful}
\end{subfigure}
\begin{subfigure}[t]{0.2\textwidth}
\includegraphics[width=\textwidth]{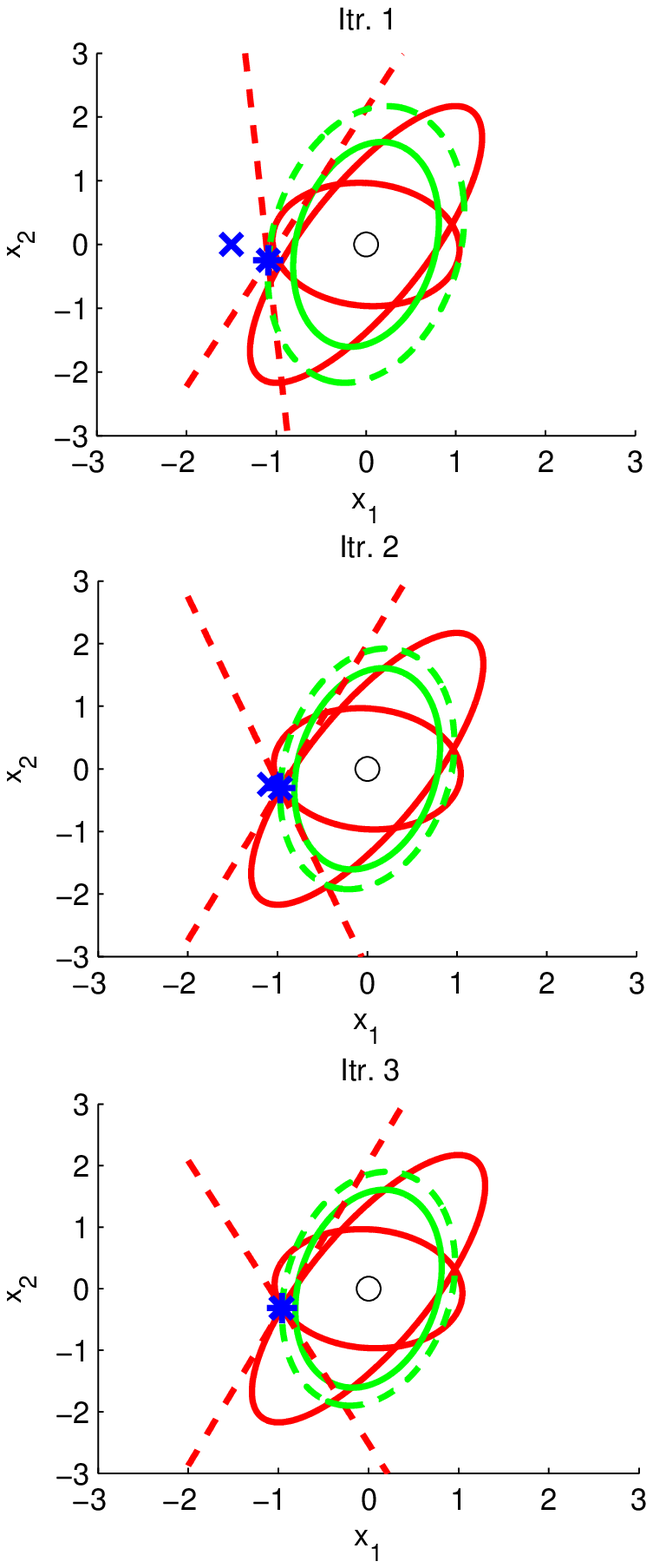}
\caption{Unsuccessful}
\end{subfigure}
\caption{Illustration of FPP-SCA algorithm in $\mathbb{R}^2$, $M=3$: 3 iterations of successful (left) and unsuccessful (right) FPP.}
\label{ExampleFig}
\end{figure}

\noindent \textbf{Illustrative example.} To get a better understanding of the approximations used in the FPP-SCA algorithm and how the solution evolves after each iteration, Fig. \ref{ExampleFig} considers a setup in $\mathbb{R}^2$ with $M=3$, where $\A_0=I$,
$\A_1 = \left[
  \begin{array}{cc}
    -1.48 & 0.68 \\
    0.68 & -0.52 \\
  \end{array}
\right]$,
 $\A_2 =  \left[
   \begin{array}{cc}
     -0.93 & -0.07 \\
     -0.07 & -1.07 \\
   \end{array}
 \right]$,
  $\A_3 =  \left[
     \begin{array}{cc}
       1.59 & -0.17 \\
       -0.17 & 0.41 \\
     \end{array}
   \right]$,
$c_1=c_2=-1$, and $c_3=1$. Note that $\A_1$ and $\A_2$ are negative semidefinite, whereas $\A_3$ is positive semidefinite. The ellipsoids that correspond to $\x^H \A_m \x = c_m$ for $m=1,2$ are plotted in red in Fig. \ref{ExampleFig}, while the ellipsoid for $m=3$ is plotted in green. 
Two different initializations for $\z_0$ are considered, both are chosen with relatively large scale so that the first two constraints of the QCQP are satisfied, but not the third one.
Therefore, only the slack variables that correspond to the third constraint can be nonzero.
The FPP-SCA algorithm is run for 3 iterations using $\lambda=10$, and the figure shows the solution after each iteration.
The point about which the non-convex constraints are linearized in each iteration is denoted by a blue cross, while the solution obtained in each iteration is denoted by a blue star.
For the non-convex constraints that correspond to $m=1,2$, the linear approximation (restriction) of each ellipsoid at each iteration is plotted in a dashed red line. 
The (extended) ellipsoid after adding the necessary slack to the convex constraint $m = 3$, in order to obtain a feasible solution for  $\large{\bf{\Pi_{3}}}$, is plotted in dashed
green. 

As shown in Fig. \ref{ExampleFig}, for the first initialization point (case (a) - left panels), a feasible solution to the original problem is obtained after the second iteration, and the optimum solution is achieved after 3 iterations. On the other hand, for the second initialization point (case (b) - right panels), the algorithm is stuck at a point that is not feasible for the original problem. Note that FPP-SCA converges in both cases, albeit to an undesirable point in the second case.

\section{Numerical Results}

To test the performance of the FPP-SCA algorithm, a problem with $n=8$ complex dimensions is considered, with $M\in\{16,24,32\}$. The entries of the matrices $\{\A_m\}_{m=1}^M$ are randomly and independently generated from a complex Gaussian distribution (with zero-mean and variance 2), then symmetrized. An initial point $\x_{\rm init}$ is randomly generated, and the values of $\{c_m\}_{m=1}^M$ are randomly generated from a Gaussian distribution $c_m \sim {\cal N}(\x_{\rm init}^H \A_m \x_{\rm init}, 1)$. If $\x_{\rm init}^H \A_m \x_{\rm init} > c_m$, the matrix $\A_m$ and $c_m$ are multiplied by $-1$ to get $\leq$ inequalities. The matrix $\A_0$ is set to the identity. To solve the SDR $\large{\bf{\Pi_{1r}}}$ and the SCA $\large{\bf{\Pi_{3}}}$, the modeling language YALMIP \cite{Yalmip} is used and  the generic conic programming solver SeDuMi \cite{Sedumi} is chosen as the solver for both approaches. The results reported in Tables \ref{SDR_Table1}-\ref{SOCP_Table2} are averaged over 1000 Monte-Carlo simulation runs.

In Table \ref{SDR_Table1}, we consider solving the QCQP $\large{\bf{\Pi_{1}}}$ using SDR $\large{\bf{\Pi_{1r}}}$ followed by a randomization (and scaling) technique that is similar to the one used in \cite{Phan09}, if the solution to  $\large{\bf{\Pi_{1r}}}$ is not rank-1. For the randomization step, $10^4$ random points are generated for each simulation run. The table reports the average number of simulation runs where a rank-1 solution was obtained, the average number of simulation runs where no feasible solution was obtained after the randomization step, the average number of simulation runs where a feasible solution was obtained with randomization, and the average difference between the solution obtained with randomization and the lower bound obtained from the (higher-rank) solution of $\large{\bf{\Pi_{1r}}}$. The table shows that as $M$ increases (i.e., the set of constraints becomes more stringent), the percentage of feasible solutions that can be obtained using the SDR approach (either directly from rank-1 solutions or after the randomization step) diminishes quickly.

\begin{table}[t!]
\centering
\caption{Results using the SDR approach for $n$ = 8.}
\begin{tabular}{|c|c | c | c|}
  \hline M	 	  					& 16 		&  24  		&   32   \\
  \hline  Rank-1 solution 	   		&  45\% 	&   16\%    &   5.5\%   \\
  \hline  No feasible sol. after randomization &  42\%		&   80\%    &   94.2\%  \\
  \hline  Feasible sol. after randomization  	&  13\%   	&    4\%    &   0.3\%  \\
  \hline   Avg. loss (dB) 			&   1.3 	&   1.5     &    -  \\
  \hline
\end{tabular}
\label{SDR_Table1}
\end{table}

\begin{table}[t!]
\centering
\caption{Results using the FPP-SCA approach for $n$ = 8.}
\begin{tabular}{|c|c|c|c|}
  \hline M	 	  				& 16 		&  24  		&   32   	\\
  \hline Feasible solution		&  100\%	&   99.5\%    &    92.8\%   \\
  \hline Avg. itrs. for feasibility&   3.207	&    4.871   &   7.893  \\
  \hline Avg. itrs. for convergence&   10.97 	&  11.6371  &   12.1352  \\
  \hline  Avg. loss (dB)	&    0.942  	&   1.5684     &    1.9256   \\
  \hline
\end{tabular}
\label{SOCP_Table1}
\end{table}

In Table \ref{SOCP_Table1}, we consider solving the QCQP $\large{\bf{\Pi_{1}}}$ using the FPP-SCA algorithm which solves  $\large{\bf{\Pi_{3}}}$ in each iteration (setting $\lambda=10$). The maximum number of iterations was set to 30 and convergence was declared if $||\x_k^H \A_0 \x_k - \x_{k-1}^H \A_0 \x_{k-1}|| \leq 10^{-4}$, for $k\geq1$.
The vector $\z_0$ used to initialize the FPP-SCA algorithm in each simulation run was randomly drawn from an i.i.d. complex circularly symmetric zero mean Gaussian distribution of variance 2.
The table reports the average number of simulation runs where a feasible solution was obtained (i.e., $\s=\mathbf{0}$), the average number of iterations until a feasible solution was obtained, the average number of iterations until convergence is declared, and the average difference between the solution obtained by FPP-SCA and the lower bound obtained from the (higher-rank) solution of $\large{\bf{\Pi_{1r}}}$. The table shows that a feasible solution can be obtained from FPP-SCA with very high probability even for large $M$, unlike the SDR approach. With $M=32$ for example, it is almost impossible to find a solution using SDR followed by randomization if the solution is not rank-1, whereas in 92.8\% of the cases FPP-SCA managed to find a feasible point. This percentage can even be increased further if multiple starting points are considered for the non-feasible cases. The table also shows that few iterations are required for the algorithm to converge, and much fewer iterations are required to reach a feasible point. Finally, the table shows that the solutions obtained using the FPP-SCA algorithm are very close to the generally unattainable relaxation lower bound provided by the SDR. Tables \ref{SDR_Table2} and \ref{SOCP_Table2} show similar results for a higher dimension $n = 20$.

\noindent {\bf Multicast Beamforming under Interference Constraints.} We further illustrate the advantage of FPP-SCA using a wireless communication design problem, namely {\em secondary multicast beamforming} as considered in \cite{Phan09}, which can be posed as the following non-convex QCQP:
\begin{equation}\label{multicast}
\begin{aligned}
\min_{{\bf w} \in \mathbb{C}^n} & \quad \|{\bf w}\|^2 \\
\text{s.t.} & \quad |{\bf w}^H{\bf h}_i|^2 \geq \tau,~~
								i = 1,...,M \\
			& \quad |{\bf w}^H{\bf g}_k|^2 \leq \eta,~~
								k = 1,...,K
\end{aligned}
\end{equation}
which describes a system comprising a secondary transmitter with $n$ antennas, $M$ secondary single-antenna receivers interested in the same multicast, and $K$ primary single-antenna receivers. The $M$ secondary receivers should be provided with signal power no less than some threshold, while the $K$ primary receivers should be protected from excessive interference. The channel gains from the transmit antennas to the $i$-th secondary user are denoted as ${\bf h}_i$, and those to the $k$-th primary user as ${\bf g}_k$. We assume i.i.d. Rayleigh fading, i.e., the channels are drawn from an i.i.d. zero-mean complex Gaussian distribution with $\sigma^2 = 1$.

For smaller problem dimensions, like the ones simulated in \cite{Phan09}, a feasible point is easy to find using SDR and randomization. However, it becomes very hard to find a feasible point when the problem size becomes higher. We conducted simulations for $n=8$, $M \in \left\{12,\cdots,24\right\}$, $K=4$, $\tau=10$ and $\eta=1$. We simulated 1000 random problem instances having feasible SDR. After SDR, we drew $10^4$ randomization points for each problem instance. None of them was (or could be scaled to be) feasible. However, FPP-SCA initialized with an SDR randomization point managed to find a feasible solution in all problem instances, with only minor average power increase (ranging from 1 to 2.2 dB, for $M$ ranging from $12$ to $24$, respectively) compared to the generally unattainable relaxation lower bound provided by SDR.


\begin{table}[t!]
\centering
\caption{Results using the SDR approach for $n$ = 20.}
\begin{tabular}{|c|c | c | c|}
  \hline M	 	  					& 32 		&  40  		&   48   \\
  \hline  Rank-1 solution 	   		&  11.5\% 	&   2.9\%    &   0.4\%   \\
  \hline  No feasible sol. after randomization &  84.7\%		&   96.4\%    &   99.5\%  \\
  \hline  Feasible sol. after randomization  	&  3.80\%   	&    0.7\%    &   0.1\%  \\
    \hline
\end{tabular}
\label{SDR_Table2}
\end{table}
\begin{table}[t!]
\centering
\vspace{1em}
\caption{Results using the FPP-SCA approach for $n$ = 20.}
\begin{tabular}{|c|c|c|c|}
  \hline M	 	  				& 32 		&  40  		&   48   	\\
  \hline Feasible solution			&  100\%	&   100\%    &    100\%   \\
  \hline Avg. itrs. for feasibility&   4.3640	&    5.0927   &  5.9840  \\
  \hline Avg. itrs. for convergence&   16.2130	&  16.5368  &   17.2805  \\
  \hline  Avg. loss (dB)	&    0.4570  	&   0.4881     &    0.5618   \\
  \hline
\end{tabular}
\label{SOCP_Table2}
\end{table}

\section{Conclusions}
FPP-SCA is a new iterative approach for approximately solving general QCQPs. 
FPP-SCA was compared with conventional SDR followed by randomization, and it was observed that FPP-SCA was successful in obtaining good feasible solutions for a much higher percentage of problem instances than SDR plus randomization, at a lower worst-case complexity, and smaller gap to the relaxation lower bound. The results suggest that FPP-SCA holds promise for a broad range of applications in engineering design problems that can be cast as non-convex QCQPs.

\end{document}